\documentclass{conm-p-l}

\usepackage{amsmath, amsthm, latexsym, amssymb, amsfonts, graphicx}

\theoremstyle{plain}
\newtheorem{Th}{Theorem}

\theoremstyle{definition}
\newtheorem{Rem}{Remark}

\newtheorem{Def}{Definition}

\newcommand{\C}{\mathbb C}
\newcommand{\Z}{\mathbb Z}
\newcommand{\R}{\mathbb R}
\newcommand{\G}{\Gamma}

\newcommand{\D}{\Delta}

\newcommand{\sig}{\sigma}

\newcommand{\res}{\operatorname{res}}
\newcommand{\vrt}{\operatorname{vert}}

\newcommand{\rs}[1]{Section~\ref{S:#1}}

\newcommand{\re}[1]{(\ref{e:#1})}

\newcommand{\rt}[1] {Theorem~\ref{T:#1}}

\copyrightinfo{2003}{(American Mathematical Society)}
\subjclass[2000]{Primary 14M25; Secondary 52B11}
\keywords{Grothendieck residues, Newton polytopes, combinatorial
coefficient}

\makeindex

\begin{document}


\title[Combinatorial coefficients, Gelfond--Khovanskii residue
formula]{On combinatorial coefficients and the Gelfond--Khovanskii
residue formula}
\author{Ivan Soprounov}
\date{}

\begin{abstract}  The Gelfond--Khovanskii residue formula computes
the sum of the values of any Laurent polynomial over solutions
of a system of Laurent polynomial equations whose Newton
polytopes have sufficiently general relative position. We discuss
two important consequences of this result: an explicit elimination
algorithm for such systems and a new formula for the mixed volume.
The integer coefficients that appear in the Gelfond--Khovanskii
residue formula are geometric invariants that depend only on
combinatorics of the polytopes. We explain how to compute
them explicitly.
\end{abstract}

\maketitle



\setcounter{page}{343}

\section{Introduction}
A Laurent polynomial is a finite linear combination of monomials
with complex coefficients:
$$f(t)=\sum_{k\in\Z^n}\lambda_kt^k,\quad t^k=t_1^{k_1}\dots t_n^{k_n},\ \
k=(k_1,\dots,k_n),\ \ \lambda_k\in\C.$$ The convex hull of those
lattice points $k\in\Z^n$ for which $\lambda_k\neq 0$ is a convex
polytope in $\R^n$, which is called the Newton polytope of $f$.
The value of $f(t)$ is defined for all $t$ in the set
$(\C\setminus 0)^n=\{t\in\C^n\ |\ t_1\dots t_n\neq 0\}$.

In this paper we talk about systems of $n$ Laurent polynomial
equations in $(\C\setminus 0)^n$ whose Newton polytopes have
sufficiently general relative positions (the precise definition
will be given in \rs{ccdef}). O.~Gelfond and A.~Khovanskii
\cite{G-Kh2} discovered that one can say a lot about solutions of
such systems explicitly in terms of the Newton polytopes and the
coefficients of the system. The following theorem (\cite{G-Kh2},
Corollary 1.6) illustrates that.

\begin{Th}\label{T:main} Consider a system of $n$ Laurent
polynomial equations
\begin{equation}\label{e:system}
f_1(t)=\dots=f_n(t)=0,\quad t\in(\C\setminus 0)^n,
\end{equation}
whose Newton polytopes $\D_1,\dots,\D_n$ have sufficiently general
relative positions. The sum of the values of a Laurent polynomial
$q$ over solutions\footnote{The assumption on the position of the
polytopes implies that the solution set of the system is finite.}
of \re{system} counted with multiplicities is equal to
$$(-1)^n\sum_{A\in\vrt(\D)}c_A\res_{A}\!\Big(q\frac{df_1}{f_1}\wedge\dots\wedge\frac{df_n}{f_n}\Big),$$
where the sum is taken over the vertices $A$ of the Minkowski sum
$\D=\D_1+\dots+\D_n$,
$\res_{A}\!\big(q\frac{df_1}{f_1}\wedge\dots\wedge\frac{df_n}{f_n}\big)$ is
the residue at the vertex $A$, and $c_A$ is the combinatorial coefficient
at $A$.
\end{Th}

The residue at a vertex $A$ is an explicit Laurent polynomial in the
coefficients of $q$ and $f_1,\dots,f_n$. In the case when $n=1$
the polytope $\D$ is a segment and the residues at the two
endpoints of $\D$ are the classical residues of $q\frac{df}{f}$ at
zero and infinity. We give the general definition in \rs{res}.

The coefficients $c_A$ are integer numbers which depend only on combinatorics of
the polytopes $\D_1,\dots,\D_n$. They can be computed by counting certain
complete flags of faces of the Minkowski sum $\D$. We discuss this
in \rs{cc}.

A particular case of \rt{main} produces a formula for the mixed
volume of the polytopes $\D_1,\dots,\D_n$ in terms of the vertices
of the Minkowski sum $\D$ and combinatorial coefficients.
Furthermore, \rt{main} provides with an explicit elimination
algorithm for systems \re{system}.

The present paper is expository. Proofs of the results that we
mention here are contained in \cite{G-Kh2,Kh,S}.

\section{Combinatorial coefficients}\label{S:cc}\index{combinatorial coefficient}

\subsection{Definition}\label{S:ccdef}

Consider a collection of $n$ convex polytopes $\D_1,\dots,\D_n$ in~$\R^n$.
Let $\D$ be their Minkowski sum,
$\D=\D_1+\dots+\D_n$. It is not hard to see that each face $\G$ of
$\D$ has a unique decomposition into the sum of faces of the $\D_i$:
\begin{equation}\label{e:decomp}
\G=\G_1+\dots+\G_n,\quad\G_i\text{ a face of }\D_i.\nonumber
\end{equation}
We will call $\G_i$ the {\it $i$-th summand} of $\G$. A face $\G$
of $\D$ is called {\it locked} if at least one of its summands has
zero dimension, i.e. is a vertex. A vertex $A$ of $\D$ is called
{\it critical} if all the proper faces of $\D$ that contain $A$
are locked.

Furthermore, consider a continuous piecewise-linear map
$\chi:\D\to\R^n$ whose $i$-th component $\chi_i$ is non-negative
and vanishes precisely on those faces $\G$ of $\D$ whose $i$-th
summand is a vertex. Such map is called {\it characteristic}. Note
that the preimage of the origin under a characteristic map $\chi$
consists of the vertices of $\D$ only. Also in a neighborhood of
each critical vertex $A$ the map $\chi$ sends the boundary
$\partial\D$ to the boundary $\partial\R^n_+$ of the positive
octant.

\begin{Def} The {\it combinatorial coefficient} $c_A$ at a critical vertex $A$
is the local degree at $A$ of the restriction of $\chi$ to the
boundary $\partial\D$:
$$\chi_A:(\partial\D,A)\to(\partial\R^n_+,0).$$
\end{Def}
One can check that it is independent of the choice of a
characteristic map $\chi$ (see \cite{Kh} for details).

The sign of the combinatorial coefficient depends on the
orientation of $\R^n$ and the order of the polytopes. Indeed, by
definition, the local degree of
$\chi_A:(\partial\D,A)\to(\partial\R^n_+,0)$ depends on the
orientation of the space containing the polytope $\D$ and on the
orientation of the target space. The latter is determined by the
order of the components $\chi_i$ of $\chi$, which corresponds to
the order of the polytopes.

Now we will give the precise definition of generic relative position
of $n$ polytopes in $\R^n$. Let $\D_1,\dots,\D_n$ be convex polytopes in $\R^n$.
Every linear functional $\xi$ on $\R^n$ defines a collection of faces
$\D_1^\xi,\dots,\D_n^\xi$ of the polytopes such that the
restriction of $\xi$ to $\D_i$ attains its minimum precisely at
$\D_i^\xi$.

\begin{Def}\label{D:generic}
We say that $n$ polytopes $\D_1,\dots,\D_n$ in $\R^n$ have {\it generic relative
position} if for any non-zero linear functional $\xi$ on $\R^n$
at least one of the faces $\D_1^\xi,\dots,\D_n^\xi$ is a vertex.
\end{Def}

Note that if $\D_1,\dots,\D_n$ have generic relative position then
every proper face of the Minkowski sum $\D$ is locked. Therefore
the combinatorial coefficient $c_A$ is defined at every vertex $A$ of $\D$.

\subsection{Explicit description}

The following theorem gives an explicit description of the
combinatorial coefficient as the number of certain complete flags
of faces of $\D$, counted with signs. The {\it sign} of a complete
flag
$$\G^0\subset\G^1\subset\dots\subset\G^{n-1}\subset\G^n,\quad\dim\G^i=i$$
of faces of a polytope in $\R^n$ is the orientation of any frame
of vectors $(v_1,\dots,v_n)$, where $v_i$ starts at $\G^0$ and
points strictly inside $\G^i$. The sign depends on the orientation
of $\R^n$.

\begin{Th}\label{T:cc}\cite{S}
Let $\D_1,\dots,\D_n$ be a collection of convex polytopes in
$\R^n$ and~$\D$ their Minkowski sum. The combinatorial
coefficient $c_A$ at a critical vertex $A$
is equal to the number of all complete flags
$$A=\G^0\subset\G^1\subset\dots\subset\G^{n-1}\subset\G^n=\D,\quad \dim\G^i=i$$
counted with signs, where $\G^i$ is an $i$-dimensional face of
$\D$ whose first $i$ summands have positive dimensions and the
last $n-i$ summands have zero dimension.
\end{Th}

Notice that this description is highly non-symmetric in the order
of the polytopes. In fact, the {\it actual} number of complete
flags can be different for different orderings, whereas the {\it
signed} number of complete flags is always the same (up to a
sign), since the combinatorial coefficient is alternating on the
order of $\D_1,\dots,\D_n$ as it follows from definition.

\subsection{Example}

Two plane polygons have generic
relative positions if and only if they do not have two parallel sides
with the same direction of an inner normal. For example,
the square and the triangle on Figure 1 have generic relative
position. Each side of the Minkowski sum comes from either one of the
polygons. The combinatorial coefficient of a vertex is 0 if
both its adjacent sides come from the same polygon, and is $\pm 1$ otherwise.
The sign depends on the orientation defined by the two adjacent sides.

\begin{figure}[ht]
\includegraphics{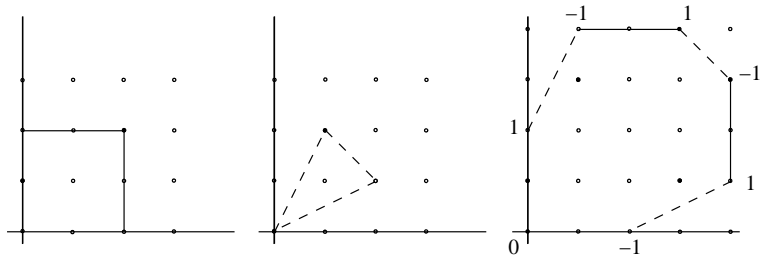}
\caption{The Minkowski sum of two polygons and combinatorial coefficients.}
\end{figure}

\begin{Rem} Combinatorial coefficients first appeared in 1996 in \cite{G-Kh}
where \rt{main} was announced. The explicit description (\rt{cc})
was found in 2001 and provided a connection between the
Gelfond--Khovanskii residue formula and Parshin's reciprocity laws
(see \cite{S}). It would be interesting to see if this description
can be useful in mixed volume computation and elimination process
which we discuss in the next section.
\end{Rem}

\section{Applications of \rt{main}}

\subsection{Mixed Volume}

The special case of \rt{main} when $q=1$ gives a formula for the
number of solutions of the system counting multiplicities. One can
prove that the residue
$\res_{A}\!\big(\frac{df_1}{f_1}\wedge\dots\wedge\frac{df_n}{f_n}\big)$
is independent of the coefficients of the polynomials $f_i$ and
equals $\det(A_1,\dots,A_n)$, where $A_i$ is the $i$-th summand of
the vertex $A$ (see \rs{ccdef}). On the other hand, according to
the BKK theorem the number of solutions is $n!$ times the mixed
volume of the Newton polytopes.\footnote{The
Bernstein--Kushnirenko--Khovanskii theorem says that the number of
solutions of a {\it generic} system with fixed Newton polytopes is
$n!$ times the mixed volume of the polytopes. However, as it was
shown by Khovanskii, if the polytopes have generic relative
position then the number of solutions (counting multiplicities) is
the same for {\it all} systems with given Newton polytopes.} This
results in a formula for the mixed volume of a collection of $n$
lattice polytopes with generic relative position:
\begin{equation}\label{e:mvol}
n!V(\D_1,\dots,\D_n)=(-1)^n\sum_{A\in\vrt(\D)}c_A\det(A_1,\dots,A_n),
\end{equation}
where the sum is taken over the vertices $A$ of the Minkowski sum
$\D$, $A_i$ is the $i$-th summand of $A$, and $c_A$ is the
combinatorial coefficient at $A$.

Let us remark that this formula holds true for non-lattice polytopes
as well. This was proved by Khovanskii \cite{Kh} using regular
polyhedral subdivisions.

\subsection{Elimination Algorithm}
According to classical elimination theory the projection of the
solution set of the system \re{system} to the $i$-th coordinate
line can be given as the solution set of a single polynomial
equation in one variable. In fact, \rt{main} provides us with an
explicit algorithm for obtaining such an equation.

For example, let us obtain a polynomial $p(t_1)$ whose roots are
the first coordinates of the solutions of the system. First,
compute the number $N$ of solutions of the system (counting
multiplicities) by taking $q=1$. This is the degree of the
polynomial $p$. Then taking $q=t_1,t_1^2,\dots,t_1^N$ we get the
sums of powers of the roots of $p$. Finally, using the Newton
formulas we can express the coefficients of the polynomial $p$ via
the sums of powers of its roots. To illustrate this procedure
we give an example in \rs{example}.

\section{Definition of $\res_{A}\!\big(
q\frac{df_1}{f_1}\wedge\dots\wedge\frac{df_n}{f_n}\big)$}\label{S:res}

Here we include the definition of the residue $\res_{A}\!\big(
q\frac{df_1}{f_1}\wedge\dots\wedge\frac{df_n}{f_n}\big)$ at a
vertex $A$ of the Minkowski sum of the Newton polytopes of
$f_1,\dots,f_n$.\index{residue!at a vertex}

Before we give the general definition let us look at the
one-dimensional case. Consider a meromorphic form $q\frac{df}{f}$,
where $q,f$ are Laurent polynomials in one variable. The possible
poles of this form are the non-zero roots of $f$, zero, and
infinity. The residue of $q\frac{df}{f}$ at a non-zero root $a$ of
$f$ is equal to $\mu(a)q(a)$, where $\mu(a)$ is the multiplicity
of the root. Since the total sum of residues is zero the sum of
the values of $q$ over the non-zero roots of $f$ counting
multiplicities is the negative sum of residues at zero and
infinity:

$$\sum_{a\neq 0, f(a)=0}\mu(a)q(a)=-\left(\res_{0}q\frac{df}{f}+\res_{\infty}q\frac{df}{f}\right).$$

Recall that $\res_{0}q\frac{df}{f}$ is the coefficient of $1/t$ in
the Laurent expansion of $q\frac{f'}{f}$ at $t=0$ and
$-\res_{\infty}q\frac{df}{f}$ is the coefficient of $1/t$ in the
Laurent expansion of $q\frac{f'}{f}$ at $t=\infty$. We denote them
by $\res_{A}q\frac{df}{f}$ and $\res_{B}q\frac{df}{f}$,
respectively, where $[A,B]$ is the Newton polytope (segment) of
$f$.

Thus we obtain  \rt{main} for $n=1$:

$$\sum_{a\neq 0, f(a)=0}\mu(a)q(a)=-\left(\res_{A}q\frac{df}{f}-\res_{B}q\frac{df}{f}\right).$$

Note that the combinatorial coefficients are $c_A=1$ and $c_B=-1$.

The following definition generalizes the notion of the Laurent
expansion at zero (infinity) to the case of several variables.

\begin{Def} \cite{G-Kh2} Let $f$ be a Laurent polynomial
$$f(t)=\sum_{k\in\Z^n}\lambda_kt^k$$ and $A$ a vertex of the
Newton polytope of $f$. Then $\lambda_A\neq 0$ and
so we can write $f=\lambda_At^A(1-h)$
for some Laurent polynomial $h$. Consider a series
\begin{equation}\label{e:series}
\frac{1}{f}=\frac{1}{\lambda_At^A(1-h)}=
\lambda_A^{-1}t^{-A}(1+h+h^2+\dots).
\end{equation}
Note that any monomial can appear only in a finite number of terms
$h^i$, thus \re{series} is a well-defined power series. Given a
Laurent polynomial $g$, the formal product of $g$ and the series
\re{series} is called the {\it Laurent expansion of $g/f$ at the
vertex $A$}.
\end{Def}

Now we are ready to define the residue $\res_{A}\!\big(
q\frac{df_1}{f_1}\wedge\dots\wedge\frac{df_n}{f_n}\big)$. Note
that the Newton polytope of the product $f_1\dots f_n$ is the
Minkowski sum $\D$ of the Newton polytopes of $f_1,\dots,f_n$.
Therefore, the Laurent expansion of $g/(f_1\dots f_n)$ is defined
at any vertex $A\in\D$ for any Laurent polynomial $g$.

\begin{Def}
The number $\res_{A}\!\big(
q\frac{df_1}{f_1}\wedge\dots\wedge\frac{df_n}{f_n}\big)$ is the
coefficient of $\frac{1}{t_1\dots t_n}$ in the Laurent expansion
of
$$\frac{qJ_f}{f_1\dots f_n},$$
at the vertex $A$. Here $J_f$ is the Jacobian
$\frac{\partial(f_1,\dots,f_n)}{\partial(t_1,\dots,t_n)}$.
\end{Def}

\section{Example}\label{S:example}

We now consider an example of a system of two polynomial equations
in two unknowns:
\begin{equation}\label{e:eqexample}
f_1(x,y)=a_{10}x+a_{01}y+a_{22}x^2y^2=0, \quad
f_2(x,y)=b_{00}+b_{12}xy^2+b_{21}x^2y=0.
\end{equation}

The following are the Newton polygons $\D_1$, $\D_2$ of $f_1$,
$f_2$ and their Minkowski sum $\D$ with the combinatorial
coefficients at the vertices.

\begin{figure}[ht]
\includegraphics{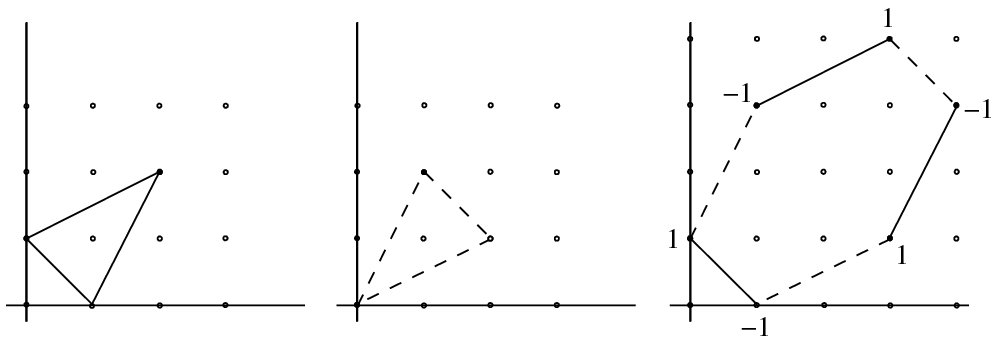}
\end{figure}

By \re{mvol} the number of solutions of the given system \re{eqexample}
(which is twice the mixed volume of $\D_1$, $\D_2$) is equal to
$$N=\left|\begin{matrix}0 & 0\\ 1 & 0\end{matrix}\right|-
    \left|\begin{matrix}0 & 1\\ 1 & 2\end{matrix}\right|+
    \left|\begin{matrix}2 & 1\\ 2 & 2\end{matrix}\right|-
    \left|\begin{matrix}2 & 2\\ 2 & 1\end{matrix}\right|+
    \left|\begin{matrix}1 & 2\\ 0 & 1\end{matrix}\right|-
    \left|\begin{matrix}1 & 0\\ 0 & 0\end{matrix}\right|=6.$$

We are going to obtain two degree six polynomials in one variable
whose roots are the $x$- and $y$-coordinates of the solution points of
\re{eqexample}, respectively.

To make our notations shorter we denote the residue
$\res_{A}\!\big( q\frac{df_1}{f_1}\wedge\frac{df_2}{f_2}\big)$ at
the vertex $A=(l,m)$ by $\res_{(l,m)}(q)$. We compute the residues
at the six vertices of $\D$ for $q=x^k$, $1\leq k\leq 6$:

\smallskip
\begin{flushleft}
${\displaystyle \res_{(0,1)}(x^k)=\res_{(1,0)}(x^k)=\res_{(1,3)}(x^k)=\res_{(3,4)}(x^k)=0,\ \ 1\leq k\leq 6}$\\
${\displaystyle \res_{(3,1)}(x^k)=\begin{cases}
    0, & k=1,2,4,5\\
    \frac{3b_{00}(a_{01}b_{21}-a_{22}b_{00})}{a_{10}b_{21}^2}, & k=3  \\
    \frac{3b_{00}^2(2a_{10}b_{12}(a_{01}b_{21}-2a_{22}b_{00})+(a_{01}b_{21}-a_{22}b_{00})^2)}{a_{10}^2b_{21}^4}, & k=6
 \end{cases}}$\\
${\displaystyle \res_{(4,3)}(x^k)=\begin{cases}
     0, & k=1,2,4,5\\
    \frac{-3b_{12}(a_{01}b_{21}-a_{10}b_{12}+2a_{22}b_{00})}{a_{22}b_{21}^2}, & k=3  \\
    \frac{-3b_{12}^2(2a_{22}b_{00}(a_{01}b_{21}-2a_{10}b_{12}+3a_{22}b_{00})+(a_{01}b_{21}-a_{10}b_{12})^2)}{a_{22}^2b_{21}^4}, & k=6
  \end{cases}}$
\end{flushleft}
\smallskip

Applying \rt{main} for $q=x^k$ we get the sums $s_k$ of the $k$-th powers of the
$x$-coordinates of the solution points of \re{eqexample}:

\smallskip
\begin{flushleft}
${\displaystyle s_1=s_2=s_4=s_5=0,}$\\
${\displaystyle s_3=\res_{(3,1)}(x^3)-\res_{(4,3)}(x^3),}$\\
${\displaystyle s_6=\res_{(3,1)}(x^6)-\res_{(4,3)}(x^6).}$
\end{flushleft}
\smallskip

We now use the Newton formulas
$$s_k-s_{k-1}\sig_1+\dots+(-1)^{k-1}s_1\sig_{k-1}+(-1)^kk\sig_k=0,\quad 1\leq k\leq 6$$
to obtain the elementary symmetric polynomials $\sig_k$:

\smallskip
\begin{flushleft}
${\displaystyle \sig_1=\sig_2=\sig_4=\sig_5=0},$\\
${\displaystyle \sig_3=\frac{s_3}{3}=
\frac{a_{01}b_{21}(a_{10}b_{12}+a_{22}b_{00})-(a_{10}b_{12}-a_{22}b_{00})^2}
{a_{10}a_{22}b_{21}^2}},$\\
${\displaystyle \sig_6=\frac{s_3^2-3s_6}{18}=
\frac{a_{01}^2b_{00}b_{12}}{a_{10}a_{22}b_{21}^2}}.$
\end{flushleft}
\smallskip

Thus we get a polynomial equation for the $x$ coordinate:

$$x^6-\frac{a_{01}b_{21}(a_{10}b_{12}+a_{22}b_{00})-(a_{10}b_{12}-a_{22}b_{00})^2}
{a_{10}a_{22}b_{21}^2}x^3+
\frac{a_{01}^2b_{00}b_{12}}{a_{10}a_{22}b_{21}^2}=0.$$

Since the system \re{eqexample} is symmetric, an equation for the
$y$ coordinate is similar:

$$y^6-\frac{a_{10}b_{12}(a_{01}b_{21}+a_{22}b_{00})-(a_{01}b_{21}-a_{22}b_{00})^2}
{a_{01}a_{22}b_{12}^2}y^3+
\frac{a_{10}^2b_{00}b_{21}}{a_{01}a_{22}b_{12}^2}=0.$$


\begin{Rem} The fact that the constant terms of the above polynomials
are so simple is not a coincidence, but a consequence of a general
result due to Khovanskii~\cite{Kh}. Namely, if the Newton polytopes
of a system have generic relative position then the product of
the $i$-th coordinates of the solution points is always a Laurent {\it monomial}
in the coefficients of the system (see \cite{Kh} for details).
\end{Rem}

\begin{Rem} Let us remark that the two equations above can be
found by computing the resultants of $f_1$ and $f_2$ with respect
to $x$ and $y$, respectively. However in higher dimensions
computing resultants can be a harder problem than finding
residues.
\end{Rem}

\begin{flushleft}
{\footnotesize{\sc Department of Mathematics and Statistics,
University of Massachusetts, Amherst, MA 01003, USA}\\
{\it E-mail address}: {\tt isoprou@math.umass.edu}}
\end{flushleft}

\end{document}